\documentclass[12pt, reqno]{amsart}
\usepackage{amssymb}

\newtheorem{theorem}{Theorem}[section]
\newtheorem{proposition}[theorem]{Proposition}
\newtheorem{lemma}[theorem]{Lemma}
\newtheorem{corollary}[theorem]{Corollary}

\theoremstyle{definition}

\newtheorem{remark}[theorem]{Remark}
\newtheorem{question}[theorem]{Question} 

\newcommand{\C}{\mathbb C}

\newcommand{\T}{\mathbb T}
\newcommand{\Z}{\mathbb Z}

\DeclareMathOperator{\pisys}{{\rm sys}\pi}
\DeclareMathOperator{\vol}{{\rm vol}} 
\DeclareMathOperator{\inj}{{\rm inj}} 

\DeclareMathOperator{\area}{{\rm area}}
\DeclareMathOperator{\length}{{\rm length}}
\DeclareMathOperator{\ent}{{\rm hom\ ent}}
\DeclareMathOperator{\card}{{\rm card}}
\DeclareMathOperator{\sigmaa}{{\rm SR}}
\DeclareMathOperator{\Emb}{{\rm Emb}}
\DeclareMathOperator{\MinEnt}{{\rm MinEnt}}
\DeclareMathOperator{\Cat}{{\rm ballcat}}

\def\gmetric{{\mathcal G}} \def\genus{{g}}

\def\ie {{\it i.e.\ }} 
 
\def\cf {\hbox{\it cf.\ }}

\numberwithin{equation}{section} 
\numberwithin{figure}{section}

\begin{document}

\title[Entropy, systoles, and asymptotic bounds]{Entropy of
systolically extremal surfaces and asymptotic
bounds$^1$}

\author[M.~Katz]{Mikhail G. Katz$^{*}$} \address{Department of
Mathematics and Statistics, Bar Ilan University, Ramat Gan 52900
Israel} \email{katzmik@math.biu.ac.il} \thanks{$^{*}$Supported by the
Israel Science Foundation (grants no.\ 620/00-10.0 and 84/03)}

\author[S.~Sabourau]{St\'ephane Sabourau} \address{Laboratoire de
Math\'ematiques et Physique Th\'eorique, Universit\'e de Tours, Parc
de Grandmont, 37400 Tours, France}

\address{Mathematics and Computer Science Department, Saint-Joseph's
University, 5600 City Avenue, Philadelphia, PA 19131, USA}

\email{sabourau@gargan.math.univ-tours.fr}


\footnotetext[1]{\large{\em Ergodic Theory and Dynamical Systems}, to
appear.  See \texttt{arXiv:math.DG/0410312}}


\begin{abstract}
We find an upper bound for the entropy of a systolically extremal
surface, in terms of its systole.  We combine the upper bound with
A.~Katok's lower bound in terms of the volume, to obtain a simpler
alternative proof of M.~Gromov's asymptotic estimate for the optimal
systolic ratio of surfaces of large genus.  Furthermore, we improve
the multiplicative constant in Gromov's theorem.  We show that every
surface of genus at least~20 is Loewner.  Finally, we relate, in
higher dimension, the isoembolic ratio to the minimal entropy.
\end{abstract}

\maketitle

\tableofcontents

\section{Entropy and systole}

We show that the volume entropy $h$ (together with A. Katok's optimal
inequality for $h$) is the ``right'' intermediary in a transparent
proof of M. Gromov's asymptotic bound for the systolic ratio of
surfaces of large genus.

In this section, we review the definitions of minimal entropy and
systole. The main results of this paper are presented in
Section~\ref{tr}.

Let~$(M,\gmetric)$ be an $n$-dimensional closed Riemannian manifold.
Denote by $(\tilde{M}, \tilde{\gmetric})$ the universal Riemannian
cover of $(M,\gmetric)$.  Choose a point~$\tilde{x}_0\in \tilde{M}$.
The volume entropy (or asymptotic volume) $h(M, \gmetric)$
of~$(M,\gmetric)$ is defined as follows:
\begin{equation} \label{eq:defent}
h(M,\gmetric) = \lim_{R \rightarrow + \infty}
\frac{\log(\vol_{\tilde{\gmetric}} B(\tilde{x}_0,R) )}{R},
\end{equation}
where $\vol_{\tilde{\gmetric}} B(\tilde{x}_0,R)$ is the volume of the
ball of radius~$R$ centered at $\tilde{x}_0\in\tilde{M}$. Since $M$ is
compact, the limit in \eqref{eq:defent} exists and does not depend on
the point $\tilde{x}_0 \in \tilde{M}$ \cite{man}.  This asymptotic
invariant describes the exponential growth rate of the volume in the
universal cover.

Define the minimal volume entropy of~$M$ as the infimum of the volume
entropy of metrics of unit volume on~$M$, or equivalently
\begin{equation}
\label{12}
\MinEnt(M) = \inf_{\gmetric} h(M,\gmetric) .  \vol(M,\gmetric) ^{\frac
{1}{n}}
\end{equation}
where $\gmetric$ runs over the space of all metrics on~$M$.

The classical result of A.~Katok ~\cite{kat} states that every
metric~$\gmetric$ on a closed surface~$M$ with negative Euler
characteristic~$\chi(M)$ satisfies the optimal inequality
\begin{equation} 
\label{eq:kat}
h(\gmetric)^2 \geq \frac{2 \pi |\chi(M)|}{\area(\gmetric)}.
\end{equation}
Inequality~\eqref{eq:kat} also holds for $\ent(\gmetric)$ \cite{kat},
as well as the topological entropy, since the volume entropy bounds
from below the topological entropy (see~\cite{man}). \\

The systole of a nonsimply connected closed Riemannian
$n$-mani\-fold~$(M,\gmetric)$ is defined as
\[
\pisys_1(M,\gmetric) = \inf_\gamma \left\{ \length(\gamma) \mid \gamma
\text{ a noncontractible loop of } M \right\}.
\]
We define the systolic ratio $\sigmaa$ of~$(M,\gmetric)$ as
\begin{equation}
\sigmaa(M,\gmetric) = \frac {\pisys_1(M,\gmetric)^n} {\vol (M,
\gmetric)},
\end{equation}
and the optimal systolic ratio of~$M$ as
\begin{equation}
\sigmaa(M) = \sup_{\gmetric} \;\sigmaa(M,\gmetric),
\end{equation}
where $\gmetric$ runs over the space of all metrics on $M$.

C.~Loewner proved the first systolic inequality.  Namely, he showed
that every metric~$\gmetric$ on the torus~$\T^2$ satisfies the
inequality
\begin{equation}
\label{17}
\pisys_1(\gmetric)^2 \leq \frac {2} {\sqrt{3}} \area(\gmetric),
\end{equation}
with equality if and only if the metric $\gmetric$ is flat, while the
group of deck transformations of $(\T^2, \gmetric)$ is a lattice
homothetic to the lattice spanned by the cube roots of unity in $\C$.
Thus, we have
\[
\sigmaa(\T^2) = \frac {2} {\sqrt{3}}.
\]
See~\cite{CK} for a recent account on systolic inequalities.
Asymptotic bounds for higher systoles are studied in \cite{Ka3}.
Systolic geometry has recently seen a period of rapid growth.  The
papers \cite{bab2, Ba3, bb} study the 1-systolic constant of
higher-dimensional manifolds.  Optimal generalisations of the Loewner
inequality are studied in~\cite{Am, BK1, BK2, BCIK2, IK, KL}.  A
generalisation of Pu's inequality is studied in~\cite{BCIK1}.  A
general framework for systolic geometry in a topological context is
proposed in \cite{KR}.  The most recent developments may be found in
\cite{Sa1, KS, KS3, Sa2, Sa3}.

In Section \ref{tr}, we present the main results of the paper.  In
Section~\ref{be}, we describe the basic estimate, based on a maximal
packing argument combined with area lower bounds for systolically
optimal surfaces.  In Section~\ref{three}, we combine the basic
estimate with Katok's inequality to prove one of our main results,
Theorem~\ref{t22}.  In Section~\ref{lo}, we prove that every surface
of genus at least 20 is Loewner.  The last section contains higher
dimensional generalisations and a proof of Theorem~\ref{t23}.

\section{The results}
\label{tr}

We will relate the minimal entropy of a closed surface to its optimal
systolic ratio.  Namely, we find an upper bound for the entropy of a
systolically extremal surface, in terms of its systole.

\begin{proposition}
Every extremal metric~$\gmetric$ on a surface~$M$ satisfies
\begin{equation} 
h(\gmetric) \leq - \frac{1}{\beta \pisys_1(\gmetric)} \log \left(
2\alpha^2 \sigmaa(M) \right),
\end{equation}
whenever $\alpha,\beta >0$ and $4\alpha + \beta < \frac{1}{2}$.
\end{proposition}

We combine this upper bound with A.~Katok's optimal lower bound in
terms of the volume, see~\eqref{eq:kat}, to obtain a simpler
alternative proof of M.~Gromov's asymptotic estimate for the optimal
systolic ratio of surfaces of large genus.  Furthermore, we improve
the multiplicative constant in Gromov's theorem (see
Section~\ref{three} for further details).

\begin{theorem}
\label{t22}
A surface~$M_\genus$ of genus~$\genus$ satisfies the bound
\begin{equation} 
\sigmaa(M_\genus) \leq \frac {\log^2 \genus}{\pi \genus} (1 + o(1))
\mbox{ when } g\to \infty .
\end{equation}
\end{theorem}

An alternative approach is taken by F. Balacheff \cite{bal}, but his
constant is not as good as Gromov's, see Section~\ref{three}.
Furthermore, the approach of the present paper lends itself to
higher-dimensional generalisations~\cite{Sa3}.

As an application, we show that every metric on a surface of genus at
least~20 satisfies the Loewner inequality~\eqref{17} for the
torus. This improves the best earlier estimate of~50.

Finally, we relate, in higher dimension, the optimal isoembolic ratio
to the minimal entropy \eqref{12}.  Recall that the optimal isoembolic
ratio of an $n$-manifold~$M$ is defined as
\begin{equation} \label{eq:defemb}
\Emb(M) = \inf_\gmetric \frac{\vol(M,\gmetric)}{\inj(M,\gmetric)^n}
\end{equation}
where $\gmetric$ runs over the space of all metrics on~$M$ and
$\inj(M,\gmetric)$ is the injectivity radius of~$(M,\gmetric)$.  We
show the following (see Section~\ref{five} for further details).

\begin{theorem}
\label{t23}
There exits a positive constant~$\lambda_n$ such that every
$n$-manifold~$M$ satisfies
\begin{equation}
\Emb(M) \geq \lambda_n \frac{\MinEnt(M)^n}{\log^n (1+\MinEnt(M))} .
\end{equation}
\end{theorem}

\section{Basic estimate}
\label{be}

In Lemma~\ref{lem:ent} we recall the following well-known fact: the
volume entropy agrees with the exponential growth rate of orbits under
the action of the fundamental group in the universal cover, sometimes
called the {\em critical exponent}.

We will also need the following estimate, \cf \eqref{eq:*}.  M.~Gromov
showed in~\cite{Gr1} that every aspherical closed surface~$M$
satisfies the inequality
\begin{equation}
\label{3/4}
\sigmaa(M) \leq \frac{4}{3}.
\end{equation}  
He also showed that every nonsimply connected closed surface admits an
extremal metric in a suitable generalized sense, namely a
metric~$\gmetric_{\rm ex}$ with optimal systolic
ratio~$\sigmaa(M,\gmetric_{\rm ex}) = \sigmaa(M)$.  Furthermore, the
disks~$D(x,r)$ of radius~$r \leq \frac{1}{2} \pisys_1(\gmetric_{\rm
ex})$ of extremal surfaces satisfy the bound
\begin{equation}
\label{1.7}
\area D(x,r) \geq 2 r^2 .
\end{equation}

We will use arguments developed in~\cite[p.~357]{kat} to prove
Proposition~\ref{prop:ent} below.  Related arguments have been
exploited in \cite{BGLM, Ge}.
\begin{proposition} 
\label{prop:ent}
Every extremal metric~$\gmetric$ on a surface~$M$ satisfies
\begin{equation} 
\label{eq:ent}
h(\gmetric) \leq - \frac{1}{\beta \pisys_1(\gmetric)} \log \left(
2\alpha^2 \sigmaa(M) \right),
\end{equation}
whenever $\alpha,\beta >0$ and $4\alpha + \beta < \frac{1}{2}$.
\end{proposition}

\begin{remark}
As shown in~\cite{Sa2}, the volume entropy of surfaces with unit
systole is bounded from above by a constant which does not depend on
the metric.  However, the constant found is not as good as in
inequality~\eqref{eq:ent}.
\end{remark}

\begin{proof}[Proof of Proposition~\ref{prop:ent}]
The idea is to bound from above the number of homotopy classes of
based loops in $M$, by deforming such a loop into a fixed subgraph in
$M$.  The subgraph can be allowed to be as ``coarse'' as the size of
the systole of $M$.  More precisely, let $x_0\in M$ be a fixed
basepoint.  Consider a maximal system of disjoint disks
\begin{equation}
\label{22}
D_i= D(x_i,R) \subset M
\end{equation}
of radius $R=\alpha \pisys_1(\gmetric)$ and centers~$x_i$ with~$i \in
I$, including $x_0$.  Since the metric~$\gmetric$ is assumed extremal,
inequality \eqref{1.7} implies
\begin{equation} 
\label{eq:*}
\area D_i \geq 2 \alpha^2 \pisys_1(\gmetric)^2 \quad \forall i\in I.
\end{equation}
Therefore, this system admits at most $\frac{\area(\gmetric)}{2 \alpha^2
\pisys_1(\gmetric)^2}$ disks.  Thus,
\begin{equation} 
\label{eq:**}
|I| \leq \left( {2 \alpha^2} {\sigmaa(M)} \right)^{-1}.
\end{equation}
Let $c:[0,T] \to M$ be a geodesic loop of length~$T$ based at~$x_0$.
Let
\begin{equation}
\label{23}
m = \left[\frac{T}{\beta \pisys_1(\gmetric)} \right]
\end{equation}
be the integer part.  The point $p_0=x_0$, together with the points
\[
p_k = c(k \beta \pisys_1(\gmetric)), \quad k=1, \ldots, m
\]
and the point $p_{m+1}=x_0$, partition the loop~$c$ into $m+1$
segments of length at most~$\beta \pisys_1(\gmetric)$.  Since the system of
disks~$D_i$ is maximal, the disks of radius $2R = 2 \alpha
\pisys_1(\gmetric)$ centered at~$x_i$ cover~$M$.  Therefore, for every~$p_k$,
a nearest point~$q_k$ among the centers~$x_i$, is at distance at
most~$2R$ from~$p_k$.  Consider the loop
\[
\alpha_k= c_k \cup [p_{k+1},q_{k+1}] \cup [q_{k+1},q_k] \cup
[q_k,p_k],
\]
where $c_k$ is the arc of~$c$ joining $p_k$ to~$p_{k+1}$, while
$[x,y]$ denotes a minimizing path joining $x$ to~$y$.  Then
\[
\length(\alpha_k) \leq 2(4 \alpha + \beta) \pisys_1(\gmetric) < \pisys_1(\gmetric),
\]
by our hypothesis on $\alpha, \beta$.  Thus the loop $\alpha_k$ is
contractible.  The same is true for the loops $c_0 \cup [p_1,q_1] \cup
[q_1,x_0]$ and $c_m \cup [x_m,q_m] \cup [q_m,p_m]$.

Therefore, the geodesic loop~$c$ is homotopic to a piecewise geodesic
loop
\begin{equation}
\label{38}
c'= (x_0,q_1, \dots, q_{m} ,x_0).
\end{equation}
Note that the minimizing path from $p_k$ to~$q_k$ may not be unique,
but we choose one, being careful that the same choice is used on
``both sides'', \ie for both loops $\alpha_{k-1}$ and $\alpha_k$.

Thus, two nonhomotopic closed geodesic loops~$c_1$ and~$c_2$ based at
$x_0$, give rise to two distinct loops $c'_1$ and $c'_2$ as in
\eqref{38}.  In constructing the loops $c'$, we always choose the same
minimizing path between a given pair of points $q_k$.  Thus, the
number $P'(T)$ of homotopy classes which can be represented by loops
of length~$T$ based at~$x_0$ satisfies
\begin{eqnarray}
P'(T) & \leq & |I|^{m} \nonumber \\ & \leq & |I|^{\frac{T}{\beta
      \pisys_1(\gmetric)}} \nonumber \\ & \leq & \left( {2 \alpha^2}
      {\sigmaa(M)} \right)^{-\frac{T}{\beta
      \pisys_1(\gmetric)}} \label{eq:***},
\end{eqnarray}
and the proposition now follows from Lemma~\ref{lem:ent}.
\end{proof}

\begin{remark}
Instead of relying on the existence of systolically extremal surfaces,
we could have exploited instead $\epsilon$-regular surfaces satisfying
\eqref{1.7} for $r\geq \epsilon$, whose existence is considerably
easier to establish~\cite[5.6.C'']{Gr1}.  Since the choice of $\alpha$ in
\eqref{3.5} entails exploiting packings by arbitrarily small disks in
\eqref{22}, we have to be careful to choose $\epsilon < \alpha
\pisys_1(M)$.
\end{remark}

\begin{remark}
A similar technique can be used to bound the number of free homotopy
classes of loops in $X$.  In higher dimensions, the latter bound seems
to be less useful, therefore we chose to bound the number of based
loops.
\end{remark}

Combining inequalities \eqref{eq:kat} and \eqref{eq:ent}, we obtain
the following corollary.

\begin{corollary}
\label{2.3}
Let $M_\genus$ be a closed orientable surface of genus $\genus$.
Whenever $4 \alpha + \beta < \frac{1}{2}$, we have
\begin{equation}
\label{35}
\frac {\log^2 \left( 2 \alpha_{\phantom{I^{i}}}^{2} \sigmaa(M_\genus)
\right)} {\sigmaa(M_\genus)} \geq 4 \pi \beta^2 (\genus -1) .
\end{equation}
\end{corollary}

We conclude this section by recalling the following well-known fact,
\cf~\cite[Proposition 9.6.6, p.~374]{KH}.

\begin{lemma} \label{lem:ent}
Let $(M,\gmetric)$ be a closed Riemannian manifold. Then,
\begin{equation}
h(M,\gmetric) = \lim_{T \rightarrow +\infty} \frac{\log(P'(T))}{T}
\end{equation}
where $P'(T)$ is the number of homotopy classes of loops based at some
fixed point~$x_0$ which can be represented by loops of length at
most~$T$.
\end{lemma}

\begin{proof}
Fix $x_0$ in $M$ and a lift $\tilde{x}_0$ in the universal
cover~$\tilde{M}$.  The group 
\[
\Gamma:=\pi_1(M,x_0)
\]
acts on~$\tilde{M}$ by isometries.  The orbit of~$\tilde{x}_0$
by~$\Gamma$ is denoted~$\Gamma.\tilde{x}_0$.  Consider a fundamental
domain~$\Delta$ for the action of~$\Gamma$, containing~$\tilde
{x}_0$. Denote by~$D$ the diameter of~$\Delta$.  The cardinal of
$\Gamma.\tilde{x}_0 \cap B(\tilde{x}_0,R)$ is bounded from above by
the number of translated fundamental domains~$\gamma.\Delta$, where
$\gamma\in \Gamma$, contained in $B(\tilde{x}_0,R+D)$.  It is also
bounded from below by the number of translated fundamental
domains~$\gamma.\Delta$ contained in~$B(\tilde{x}_0,R)$.  Therefore,
we have
\begin{equation}
\frac{\vol(B(\tilde{x}_0,R))}{\vol(M,\gmetric)} \leq \card
(\Gamma.\tilde{x}_0 \cap B(\tilde{x}_0,R)) \leq
\frac{\vol(B(\tilde{x}_0,R+D))}{\vol(M,\gmetric)}.
\end{equation}
Take the log of these terms and divide by ${R}$.  The lower bound
becomes
\begin{equation}
\begin{aligned}
\frac{1}{R} & \log \left( \frac{\vol(B(R))}{\vol(\gmetric)} \right) =
\cr & = \frac{1}{R} \log \left( \vol(B(R)) \right) - \frac{1}{R} \log
\left( \vol(\gmetric) \right),
\end{aligned}
\end{equation}
and the upper bound becomes
\begin{equation}
\begin{aligned}
\frac{1}{R} & \log \left( \frac{\vol(B(R+D))}{\vol(\gmetric)} \right)
= \cr & = \frac{R+D}{R} \frac{1}{R+D} \log(\vol(B(R+D))) - \frac{1}{R}
\log(\vol(\gmetric)).
\end{aligned}
\end{equation}
Hence both bounds tend to~$h(\gmetric)$ when $R$ goes to infinity.
Therefore,
\begin{equation}
h(\gmetric) = \lim_{R \rightarrow +\infty} \frac{1}{R} \log \left( \card
(\Gamma .\tilde{x}_0 \cap B(\tilde{x}_0,R))\right).
\end{equation}
This yields the result since $P'(R) = \card (\Gamma.\tilde{x}_0 \cap
B(\tilde{x}_0,R))$.
\end{proof}

\section{Asymptotic behavior of systolic ratio for large genus}
\label{three}

We now consider the asymptotic behavior of the optimal systolic ratio
of surfaces.  M.~Gromov \cite[p.~74]{Gr1} established a bound for the
ratio~$\sigmaa(M_\genus)$ by using a technique known as ``diffusion of
chains''.  The multiplicative constant $\frac{1}{\pi}$ in~\eqref
{eq:asymp} below improves the constant~$4$ which could be obtained
from the techniques in~\cite{Gr1}.  F.~Balacheff~\cite{bal} found
another proof of a similar inequality, by combining the works of
S.~Kodani~\cite{Ko}, and B.~Bollob\'as and E.~Szemer\'edi~\cite{BoS}
on systolic inequalities of graphs, obtaining a multiplicative
constant of~$\frac{8}{3 (\log 2)^2}$.
\begin{theorem} 
\label{theo:asymp}
Given a real number $\lambda < \pi$, every surface~$M_\genus$ of
genus~$\genus$ satisfies
\begin{equation}
\label{31}
\sigmaa(M_\genus)  \leq \frac { \log^2 \genus} {\lambda \genus}
\end{equation}
for $\genus$ large enough.  Thus,
\begin{equation} 
\label{eq:asymp}
\sigmaa(M_\genus) \lesssim \frac {\log^2 \genus}{\pi \genus}\quad {\rm
as} \; \genus \to \infty.
\end{equation}
\end{theorem}

\begin{remark}
A lower bound for~$\sigmaa(M_{\genus})$ was found by P.~Buser and
P.~Sarnak~\cite{BS}.  Namely, they construct hyperbolic
surfaces~$(M_\genus , \gmetric_\genus)$ of arbitrarily large
genus~$\genus$ obtained as congruence coverings of an arithmetic
Riemann surface, such that
\begin{equation}
\sigmaa(M_\genus , \gmetric_\genus) \gtrsim \frac {4} {9 \pi}
\frac{\log^2 \genus} {\genus}.
\end{equation}
Therefore, we have
\begin{equation}
\frac{4}{9 \pi} \le \limsup_{\genus\to\infty} \, \sigmaa(M_\genus)
\frac{\genus} {\log^2 \genus} \le \frac{1} {\pi} .
\end{equation}
\end{remark}

\begin{proof}[Proof of Theorem~\ref{theo:asymp}]
Let $\gmetric$ be an extremal metric on~$M_\genus$.  We now apply
Corollary~\ref{2.3}.  Note that $\sigmaa(M_\genus, \gmetric)$ tends to
zero as the genus~$\genus$ becomes unbounded, \cf \eqref{1.1}.
For a right choice of $\alpha$ and $\beta$, inequality~\eqref{35}
leads to the asymptotic implicit upper bound~\eqref{36} on~$\sigmaa(M_g)$ below.
Indeed, given $\lambda < \lambda_+ < \pi$, we set
\[
\beta= \sqrt{\frac{\lambda_+}{ 4\pi}} < \frac{1}{2},
\]
and choose 
\begin{equation}
\label{3.5}
\alpha < \frac{1}{4} \left( \frac{1}{2} - \beta \right).
\end{equation}
Inequality~\eqref {35} implies
\begin{equation}
\label{36}
\frac{\log^2 (\sigmaa(M_\genus)) } {\sigmaa(M_\genus) } \geq \lambda \genus
\end{equation}
if $\genus$ is large enough.  Now we want to invert this relation in
order to get an asymptotic upper bound on~$\sigmaa(M_g)$. Let $\rho=
\sigmaa(M_\genus)^{-\frac{1}{2}}$ and $\delta=
\frac{1}{2}\sqrt{\lambda\genus}$.  Then inequality~\eqref{36} yields
the following estimate:
\[
\begin{aligned}
\rho \log \rho & \geq \delta \\ & \geq \delta - \frac{\delta \log \log
\delta}{\log \delta} \\ & = \frac{\delta}{\log \delta} \log \left(
\frac{\delta}{\log \delta} \right).
\end{aligned}
\]
Since the function $x \log x$ is increasing for $x$ large enough, we
deduce that $\rho \geq \frac{\delta}{\log\delta}$, and the latter
inequality translates back into \eqref{31}.
\end{proof}

\section{When is a surface Loewner?}
\label{lo}

We now extend the classical Loewner inequality~\eqref{17} on the torus
to surfaces of higher genus.  We will say that a surface~$M$ is {\em
Loewner\/} if $\sigmaa(M) \leq \frac{2}{\sqrt{3}}$.

\begin{theorem}
\label{41}
Every surface of genus at least~20 is Loewner.
\end{theorem}

\begin{proof}
Let $\gmetric$ be an extremal metric on~$M$.  If $\beta = \frac{1}{2}
- 4 \alpha$, inequality~\eqref{35} yields
\begin{equation} \label{eq:1}
\frac{\log^2 \left( 2 \alpha^2 \sigmaa(M) \right)} {\sigmaa(M)} \geq
4 \pi (\frac{1}{2} - 4 \alpha)^2 (\genus -1)
\end{equation}
for every $\alpha \leq \frac{1}{8}$.

Suppose now that $\sigmaa(M) > \frac{2}{\sqrt{3}}$.  Since $\sigmaa(M)
\leq \frac{4}{3}$ by \eqref{3/4}, we have $2 \alpha^2 \sigmaa(M) \leq
1$ for every $\alpha \leq \frac{1}{8}$. Therefore,
\begin{equation}
\frac{\sqrt{3}}{2} \log^2(\frac{\sqrt{3}}{4 \alpha^2}) \geq 4 \pi
(\frac{1}{2} - 4 \alpha)^2 (\genus -1)
\end{equation}
for every $\alpha \leq \frac{1}{8}$.
Hence,
\begin{equation}
\label{43}
\min_{0 < \alpha \leq \frac{1}{8}} \frac{\sqrt{3}}{8 \pi} \left(
\frac{\log \left( \frac{\sqrt{3}}{4 \alpha^2} \right)}{\frac{1}{2} - 4
\alpha} \right)^2 \geq \genus -1
\end{equation}
For $\alpha = .031$, the expression to minimize is about~$18.201$.
Thus, $\genus \leq 19$.
Therefore, every surface of genus greater or equal to~20 is Loewner.
\end{proof}

Note that M.~Gromov~\cite[p.~50]{Gr1} (\cf \cite[Theorem 4, part
(1)]{Ko}) proved a general estimate which implies that
\begin{equation}
\label{1.1}
\sigmaa(M_\genus) < \frac {64} {4\sqrt{\genus} + 27} .
\end{equation}

It follows from Gromov's estimate \eqref{1.1} that orientable surfaces
satisfy Loewner inequality \eqref{17} if the genus $\genus$ is bigger
than 50.  Our theorem brings this bound down to 20.  It has recently
been shown~\cite{KS} that the genus~2 surface is Loewner.  The
remaining open cases are therefore $\genus= 3, \ldots, 19$.

\begin{remark}
Let $\alpha = \frac{1}{30}$.  Instead of taking an arbitrary packing,
we start with a systolic loop, and choose 15 disjoint disks centered
at equally spaced points of the loop.  By an averaging procedure
\cite{Gr1}, we can get the combined area of these disks to be, not $15
(2r^2)$, but rather~$15 (3r^2)$.  Now we complete these 15 disks to a
maximal packing, and argue as before.  The only difference is that we
have a better lower bound for the area of the disks.

If $\alpha= \frac{1}{30}$, then there are at most 382 balls in the
packing of an unloewner surface.  Hence the area of the packing is at
least $2.039 \alpha^2$, instead of $2 \alpha^2$.  Calculating the
resulting expression in \eqref{43} yields about 18.12, which is better
than 18.20 but not enough to dip under~18.
\end{remark}

\section{Isoembolic ratio, minimal entropy, and simplicial norm}
\label{five}

Analogous estimates can be proved in higher dimension.  However, the
results we obtain are weaker than in the previous sections, in the
absence of similar results on the existence of systolically extremal
\mbox{metrics}.



M.~Berger proved in~\cite{ber} that the isoembolic ratio~$\Emb(M,
\gmetric)$ of every Riemannian $n$-manifold~$(M,\gmetric)$, defined as
\begin{equation}
\Emb(M,\gmetric) = \frac{\vol(M,\gmetric)}{\inj(M,\gmetric)^n},
\end{equation}
satisfies
\begin{equation}
\Emb(M,\gmetric) \geq \Emb(S^n,can)
\end{equation}
with equality if and only if $M$ is a sphere, while $\gmetric$ has
constant sectional curvature. In particular, $\Emb(M) \geq \Emb(S^n)$,
where $\Emb$ is the optimal isoembolic ratio defined
in~\eqref{eq:defemb}.

C.~Croke showed in~\cite{cro} that for every Riemannian
manifold~$(M,\gmetric)$,
\begin{equation}
\label{eq:cro}
\vol B(R) \geq c_n R^n
\end{equation}
for every $R \leq \frac{1}{2} \inj(\gmetric)$, where $\vol B(R)$ is the
volume of a ball of radius~$R$ in~$M$.
Furthermore, he showed in~\cite{cro88} that
\begin{equation} \label{eq:cat}
\Emb(M) \geq C_n \Cat(M) ,
\end{equation}
where $\Cat(M)$ is the minimum number~$n$ such that $M$ can be covered
by~$n+1$ topological balls. If $M$ is a closed $n$-dimensional
manifolds, then $\Cat(M) \leq n$ (see~\cite[p.~77]{CLOT}).  Note also
that T.~Yamaguchi~\cite{ya88} proved that, for every real~$C$, the
class of closed $n$-dimensional manifolds~$M$ satisfying $\Emb(M) \leq
C$ contains finitely many homotopy types. Furthermore,
I.~Babenko~\cite[Theorem 1.2, p.~6]{bab} showed that the minimal
volume entropy is a homotopy invariant of~$M$.\\

The following analogue of Proposition~\ref{prop:ent} holds.

\begin{proposition} \label{prop:inj}
Every Riemannian $n$-manifold~$(M,\gmetric)$ satisfies
\begin{equation} \label{eq:inj}
h(\gmetric) \leq \frac{1}{\beta \inj(\gmetric)} \log \left(
\frac{\Emb(\gmetric)}{c_n \alpha^n} \right)
\end{equation}
whenever $\alpha,\beta >0$ and $4\alpha + \beta < \frac{1}{2}$.
\end{proposition}

The proof is identical to that of Proposition~\ref{prop:ent}, with
$c_n$ from \eqref{eq:cro} replacing the coefficient $2$ in
\eqref{eq:*}, \eqref{eq:**} and~\eqref{eq:***}, and $n$ replacing the
dimension~2.

Arguing as in Theorem~\ref{theo:asymp}, we obtain the following
theorem (compare with~\eqref{eq:cat}).

\begin{theorem}
There exits a positive constant~$\lambda_n$ such that every
$n$-manifold~$M$ satisfies
\begin{equation} \label{eq:emb}
\Emb(M) \geq \lambda_n \frac{\MinEnt(M)^n}{\log^n (1+\MinEnt(M))}
\end{equation}
\end{theorem}

G.~Besson, G.~Courtois, and S.~Gallot~\cite{bcg} proved that the
minimal entropy of a closed negatively curved locally symmetric
$n$-manifold~$(M,\gmetric_0)$ satisfies
\begin{equation} 
\label{eq:bcg}
\MinEnt(M)^n = \vol(\gmetric_0).h(\gmetric_0)^n.
\end{equation}

This theorem generalizes \eqref{eq:kat}.

\begin{remark}
The relation~\eqref{eq:bcg} is a sharp version, for
negatively curved locally symmetric manifolds, of the following
general result of M.~Gromov (see~\cite[p.~37]{gro81}).  Every
$n$-manifold~$M$ with simplicial volume~$||M||$ satisfies
\begin{equation} \label{eq:gro}
\MinEnt(M)^n \geq C_n ||M||.
\end{equation}
Therefore, inequalities~\eqref{eq:emb} and~\eqref{eq:gro} show that
there exists a positive constant~$\lambda'_n$ such that
\begin{equation} \label{eq:emb/simpl}
\Emb(M) \geq \lambda'_n \frac{||M||}{\log^n (1+||M||)}.
\end{equation}
M.~Gromov proved in~\cite[p.~74]{Gr1} that there exists a positive
constant~$\lambda_n''$ such that
\begin{equation} \label{eq:simpl}
\sigmaa(M) \leq \lambda''_n \frac{\log^n (1+||M||)}{||M||}.
\end{equation}
Thus, inequality~\eqref{eq:emb/simpl} appears as a particular case
of~\eqref{eq:simpl}.
\end{remark}

\begin{question}
Are there manifolds with large minimal entropy and small simplicial
volume?  Such manifolds (if they exist) would provide examples where
\eqref{eq:emb} yields a better estimate than~\eqref{eq:simpl}.
\end{question}

More generally, in the presence of a lower bound for volumes of balls,
our argument yields the following result.

\begin{proposition}
Let~$\gmetric$ be a metric on a closed $n$-manifold~$M$ such that for some
$c>0$, every ball of radius~$r$ with $0 < r < \frac{1}{2}
\pisys_1(\gmetric)$ satisfies
\begin{equation}
\area B(r) \geq c r^n .
\end{equation}
Then, we have
\begin{equation}
h(\gmetric) \leq - \frac{1}{\beta \pisys_1(\gmetric)} \log \left( c
\alpha^n \sigmaa(M,\gmetric) \right),
\end{equation}
whenever $\alpha,\beta >0$ and $4\alpha + \beta < \frac{1}{2}$.
Therefore, there exists a positive constant~$\lambda =\lambda(n,c)$
such that
\begin{equation}
\sigmaa(M,\gmetric) \geq \lambda(n,c) \frac{\MinEnt(M)^n}{\log^n
(1+\MinEnt(M))}.
\end{equation}
\end{proposition}


\end{document}